\documentclass{amsart}
\usepackage{wasowicz}

\info{%
 Published: Studia Math. 120 (1996), 75--82.
 A manuscript was originally written in January 1996.\\
 This is a version rewritten in 2008 using the \texttt{amsart} document class.}

\begin{document}

\title{Polynomial selections and separation by polynomials}

\author{\SW}
\address{\SWaddr}
\email{\SWmail}


\date{January 1996}

\begin{abstract}
 K. Nikodem and the present author proved in \cite{NikWas95} a
 theorem concerning separation by affine functions. Our purpose is to
 generalize that result for polynomials. As a consequence we obtain two
 theorems connected with separation of $n$-convex function from $n$-concave
 function by a polynomial of degree at most $n$ and a stability result of
 Hyers-Ulam type for polynomials.
\end{abstract}

\maketitle

\section{Introduction}

By $\R$, $\N$ we denote the set of all reals and positive integers,
respectively. Let $I\subset\R$ be an interval. In this paper we present a
necessary and sufficient condition under which two functions
$f,g:I\to\R$ can be separated by a polynomial of degree at most $n$, where
$n\in\N$ is a fixed number. Our main result is a generalization of the
theorem concerning separation by affine functions obtained recently by
K. Nikodem and the present author in \cite{NikWas95}. To get it we use
Behrends and Nikodem's abstract selection theorem
(cf.~\cite[Theorem~1]{BerNik95}). It is a variation of Helly's theorem
(cf.~\cite[Theorem 6.1]{Val64}).

By $\ccR$ we denote the family of all non-empty compact real intervals.
Recall that if $F:I\to \ccR$ is a set-valued function then a function
$f:I\to\R$ is called a \emph{selection} of $F$ iff $f(x)\in F(x)$ for every
$x\in I$.

Behrends and Nikodem's theorem states that if $\mathscr{W}$ is an
$n$-dimensional space of functions mapping $I$ into $\R$ then a set-valued
function $F:I\to \ccR$ has a selection belonging to $\mathscr{W}$ if and only
if for every $n+1$ points $x_{1},\dots,x_{n+1}\in I$ there exists a
function  $f\in\mathscr{W}$ such that $f(x_{i})\in F(x_{i})$ for
$i=1,\dots,n+1$.

Let us start with the notation used in this paper. Let $n\in\N$. If
$x_{1},\dots,x_{n}\in I$ are different points then for $i=1,\dots,n$ we
define
\[c_{i}(x;x_{1},\dots,x_{n})=
  \prod_{\substack{j=1\\ j\ne i}}^{n}
  \frac{x-x_{j}}{x_{i}-x_{j}}.\]

Note that $c_{i}(x_{j};x_{1},\dots,x_{n})$ is equal to 0 if $i\ne j$ and
to 1 if $i=j$, $i,j=1,\dots,n$. ${\mathscr{P}}_{n}$ stands for
the family of all polynomials of degree at most $n$. If $x_{1},\dots,
x_{n+1}\in I$ are different points then the Lagrange interpolating
polynomial going through the points $(x_{i},y_{i})$, $i=1,\dots,n+1$, is
given by the following formula:
\begin{equation}\label{Lagrange}
 w(x)=\sum_{i=1}^{n+1} c_{i}(x;x_{1},\dots,x_{n+1})y_{i}.
\end{equation}
This uniquely determined polynomial is a member of ${\mathscr{P}}_{n}$.
Moreover, if $x<x_{1}<\dots<x_{n+1}$ then $c_{i}(x;x_{1},\dots,x_{n+1})$
is positive if $i$ is odd and negative if $i$ is even.

\section{Polynomial selections of set-valued functions}

Now we shall prove a selection theorem which will be used to obtain our main
result. If $n\in\N$ and $A_{i}\subset\R$, $i=1,\dots,n$, then
\[\sum_{i=1}^{n} A_{i}\] denotes the algebraic sum of the sets
$A_{i}$, $i=1,\dots,n$.
\begin{thm}\label{tw1}
 Let $n\in\N$. The set-valued function $F:I\to \ccR$ has a selection
 belonging to ${\mathscr{P}}_{n}$ if and only if for every $x_{0},x_{1},\dots,
 x_{n+1}\in I$ such that $x_{0}<x_{1}<\dots<x_{n+1}$ the following
 condition holds:
 \begin{equation}\label{selekcja}
  F(x_{0})\cap \biggl(\sum_{i=1}^{n+1} c_{i}(x_{0};x_{1},\dots,x_{n+1})
  F(x_{i})\biggr)\ne\emptyset.
 \end{equation}
\end{thm}

\begin{proof}
If $F$ has a selection belonging to ${\mathscr{W}}_{n}$ then
\eqref{selekcja} is obvious. We prove that \eqref{selekcja} implies an
existence of a polynomial selection of $F$. First we note that
${\mathscr{P}}_{n}$ is an $(n+1)$-dimensional space of functions. If we prove
that for every $n+2$ points $x_{0},x_{1},\dots,x_{n+1}\in I$ there exists
a $w\in {\mathscr{P}}_{n}$ such that $w(x_{i})\in F(x_{i})$, $i=0,1,\dots,n+1$,
then by Behrends and Nikodem's theorem $F$ will have a desired selection.
(For another Helly-type theorem which may be used here cf. also
\cite[Theorem 6.9]{Val64}.)

Fix any different points
 $x_{0},x_{1},\dots,x_{n+1}\in I$
and assume $x_{0}< x_{1}<\dots<x_{n+1}$.
Let $L_{i}=c_{i}(x_{0};x_{1},\dots,x_{n+1})$,
$i=1,\dots,n+1$. Thus \eqref{selekcja} has the form
\begin{equation}\label{sel1}
 F(x_{0})\cap\biggl(\sum_{i=1}^{n+1} L_{i}F(x_{i})\biggr)\ne\emptyset.
\end{equation}
It is easy to observe that $L_{i}$ is positive if $i$ is odd and negative
if $i$ is even.

Put \[y_{0}=\inf F(x_{0}),\;\;\;z_{0}=\sup F(x_{0})\]
and for $i=1,\dots,n+1$,
\[
 y_{i}=\begin{cases}
        \inf F(x_{i})&\text{if }L_{i}>0,\\
        \sup F(x_{i})&\text{if }L_{i}<0,
       \end{cases}\qquad
 z_{i}=\begin{cases}
        \sup F(x_{i})&\text{if }L_{i}>0,\\
        \inf F(x_{i})&\text{if }L_{i}<0.
       \end{cases}
\]

Therefore $F(x_{0})=\left[ y_{0},z_{0}\right]$ and for $i=1,\dots,n+1$,
\[
 F(x_{i})=\begin{cases}
           [y_{i},z_{i}]&\text{if }L_{i}>0,\\
           [z_{i},y_{i}]&\text{if }L_{i}<0.
          \end{cases}
\]
Since $-\left[\alpha,\beta\right] =\left[-\beta,-\alpha\right]$ for all
$\alpha,\beta\in\R$, we have $L_{i}F(x_{i})=\left[ L_{i}y_{i}, L_{i}z_{i}
\right]$, $i=1,\dots,n+1$. If $u=L_{1}y_{1}+\dots+L_{n+1}y_{n+1}$ and
$v=L_{1}z_{1}+\dots+L_{n+1}z_{n+1}$ then $u\le v$. Furthermore,

\begin{multline*}
 \sum_{i=1}^{n+1} L_{i}F(x_{i})
  =[L_{1}y_{1}, L_{1}z_{1}]+\dots+[L_{n+1}y_{n+1},L_{n+1}z_{n+1}]\\
  =[ L_{1}y_{1}+\dots+L_{n+1}y_{n+1},L_{1}z_{1}+\dots+L_{n+1}z_{n+1}]=[u,v]
\end{multline*}
and by \eqref{sel1} we get
\begin{equation}\label{sel2}
 \left[ y_{0},z_{0}\right]\cap\left[ u,v\right]\ne\emptyset.
\end{equation}

There are three cases of the location of the above intervals:
\begin{enumerate}[(a)]
 \item $u\in[y_{0},z_{0}]$,
 \item $v\in[y_{0},z_{0}]$,
 \item $[y_{0},z_{0}]\subset[u,v]$.
\end{enumerate}

Fix $t\in\left[ 0,1\right]$ and consider the polynomial
$\varphi _{t}\in{\mathscr{P}}_{n}$ going through $n+1$ different points:
\[
 \bigl(x_{0},tu+(1-t)v\bigr)\text{\;\;and\;\;}
 \bigl(x_{i},ty_{i}+(1-t)z_{i}\bigr)\text{\;\;for\;\;}i=1,\dots,n-1,n+1.
\]
We shall show later that
\begin{equation}\label{fi_t}
 \varphi _{t}(x_{n})=ty_{n}+(1-t)z_{n}.
\end{equation}

Hence, in the case (a) for $w=\varphi _{1}$ we have
\begin{align*}
 w(x_{0})&=u\in\left[ y_{0},z_{0}\right] =F(x_{0}),\\
 w(x_{i})&=y_{i}\in F(x_{i}),\;\;i=1,\dots,n-1,n,n+1
\end{align*}
and similarly in the case (b) for $w=\varphi _{0}$. In the case
(c) $y_{0}=\lambda u+(1-\lambda)v$ for some $\lambda\in\left[
0,1\right]$. For $w=\varphi _{\lambda}$ we obtain
\begin{align*}
 w(x_{0})&=y_{0}\in F(x_{0}),\\
 w(x_{i})&=\lambda y_{i}+(1-\lambda)z_{i}\in F(x_{i}),\;\;i=1,\dots,n-1,n,n+1.
\end{align*}
So in all cases there exists a $w\in {\mathscr{P}}_{n}$ such that $w(x_{i})\in
F(x_{i})$, $i=0,\dots,n+1$, which will complete the proof if we show that
\eqref{fi_t} holds true.

By \eqref{Lagrange} we get
\begin{align*}
 \varphi_{t}(x)&=c_{0}(x;x_{0},x_{1},\dots,x_{n-1},x_{n+1})\bigl(tu+(1-t)v\bigr)\\
               &+\sum_{i=1}^{n-1}c_{i}(x;x_{0},x_{1},\dots,x_{n-1},x_{n+1})\bigl(ty_{i}+(1-t)z_{i}\bigr)\\
               &+c_{n+1}(x;x_{0},x_{1},\dots,x_{n-1},x_{n+1})\bigl(ty_{n+1}+(1-t)z_{n+1}\bigr).
\end{align*}
If $M_{i}=c_{i}(x_{n};x_{0},x_{1},\dots,x_{n-1},x_{n+1})$, $i=0,1,\dots,
n-1,n+1$, then after a bit of computation
\[\varphi _{t}(x_{n})=\sum_{\stackrel{\scriptstyle i=1}{i\ne n}}^{n+1}
  (M_{0}L_{i}+M_{i})\bigl(ty_{i}+(1-t)z_{i}\bigr)+M_{0}L_{n}\bigl(ty_{n}+(1-t)z_{n}\bigr).\]

One can verify (using the product formula given in Introduction) that
$M_{0}L_{n}=1$ and $M_{0}L_{i}+M_{i}=0$ for $i=1,\dots,n-1,n+1$. Hence
\eqref{fi_t} holds and this finishes the proof.
\end{proof}

As a consequence of Theorem~\ref{tw1} we obtain
\begin{cor}\label{wn1}{\rm \cite[Theorem 1]{Was95}}
 A set-valued function $F:I\to \ccR$ has an affine selection iff for
 every $x,y\in I$, $t\in \left[ 0,1\right]$
 \[ F\bigl(tx+(1-t)y\bigr)\cap\bigl(tF(x)+(1-t)F(y)\bigr)\ne\emptyset.\]
\end{cor}

\begin{proof}
The above condition is equivalent to \eqref{selekcja} for
$n=1$, $x<y$, $x_{0}=x$, $x_{2}=y$, $x_{1}=tx_{0}+(1-t)x_{2}$, where
$t=\frac{x_{1}-x_{2}}{x_{0}-x_{2}}$.
\end{proof}

\section{Separation by polynomials}

The main result of this paper reads as follows
\begin{thm}\label{tw2}
 Let $n\in\N$, $f,g:I\to\R$. The following conditions are equivalent:
 \begin{enumerate}[\upshape (i)]
  \item\label{tw2:i}
   there exists $w\in {\mathscr{P}}_{n}$ such that $f(x)\le w(x)\le g(x)$, $x\in I$;
  \item\label{tw2:ii}
   $f(b)\le g(b)$, where $b\in I$ is the right side endpoint of $I$
   (if exists) and for every $x_{0},x_{1},\dots,x_{n+1}\in I$
   such that $x_{0}\le x_{1}<\dots<x_{n+1}$
 \end{enumerate}
 \begin{equation}\label{oddziel}
  \begin{aligned}
    f(x_{0})&\le \sum_{\substack{i=1\\ i\text{\ \normalfont{odd}}}}^{n+1}
       c_{i}(x_{0};x_{1},\dots,x_{n+1})g(x_{i})+
       \sum_{\substack{i=1\\ i\text{\ \normalfont{even}}}}^{n+1}
       c_{i}(x_{0};x_{1},\dots,x_{n+1})f(x_{i}),\\
     g(x_{0})&\ge \sum_{\substack{i=1\\ i\text{\ \normalfont{odd}}}}^{n+1}
       c_{i}(x_{0};x_{1},\dots,x_{n+1})f(x_{i})+
       \sum_{\substack{i=1\\ i\text{\ \normalfont{even}}}}^{n+1}
       c_{i}(x_{0};x_{1},\dots,x_{n+1})g(x_{i}).
  \end{aligned}
 \end{equation}
\end{thm}

\begin{proof}
To prove that \itemref{tw2:i} implies \itemref{tw2:ii} fix any
$x_{0},x_{1},\dots,x_{n+1}\in I$ such that
$x_{0}\le x_{1}<\dots<x_{n+1}$.
Since the polynomial $w$ goes through the points
$(x_{i},w(x_{i}))$, $i=1,\dots,n+1$, we have
\[ w(x_{0})=\sum_{i=1}^{n+1} c_{i}(x_{0};x_{1},\dots,x_{n+1})w(x_{i}).\]
Then the ineualities \eqref{oddziel} are obvious.

To prove the converse implication first note that replacing $x_{0}$ by
$x_{1}$ in \eqref{oddziel} we have $f(x_{1})\le g(x_{1})$ in both
ineualities, i.e. (ii) yields $f\le g$ on $I$. Let
\[F(x)=\left[ f(x),g(x)\right],\;\;\;x\in I.\]
We shall show that $F:I\to \ccR$ fulfils \eqref{selekcja}. Fix any $x_{0},
x_{1},\dots,x_{n+1}\in I$ such that $x_{0}<x_{1}<\dots<x_{n+1}$. Let
$u,v$ be equal to the right hand sides of the upper and lower
ineualities \eqref{oddziel}, respectively. Therefore $v\le u$ and
\begin{equation}\label{oddz1}
 \left[ f(x_{0}),g(x_{0})\right]\cap\left[ v,u\right]\ne\emptyset
\end{equation}
(otherwise $g(x_{0})<v$ or $u<f(x_{0})$ -- contradiction with
\eqref{oddziel}). Let $L_{i}=c_{i}(x_{0};x_{1},\dots,x_{n+1})$,
$i=1,\dots,n+1$. Then
\[
 L_{i}F(x_{i})=\begin{cases}
                [L_{i}f(x_{i}),L_{i}g(x_{i})]&\text{if $i$ is odd},\\
                [L_{i}g(x_{i}),L_{i}f(x_{i})]&\text{if $i$ is even}
               \end{cases}
\]
and
\[\left[ v,u\right] =\sum_{i=1}^{n+1} L_{i}F(x_{i}).\]
Thus \eqref{oddz1} implies \eqref{selekcja}. By Theorem~\ref{tw1} $F$
has a selection $w\in {\mathscr{P}}_{n}$. This finishes the proof.
\end{proof}

\begin{rem}
 Ineualities \eqref{oddziel} in Theorem~\ref{tw2} do not guarantee
 $f(b)\le g(b)$, where $b\in I$ is the right side endpoint of $I$
 (if exists). The following two functions
 \[
  f(x)=\begin{cases}
         \frac{1}{2}x&\text{for }0\le x<1,\\
         1&\text{for }x=1
        \end{cases}\qquad\text{and}\qquad
  g(x)=\begin{cases}
        x&\text{for }0\le x<1,\\
        \frac{1}{2}&\text{for }x=1
       \end{cases}
  \]
 fulfil \eqref{oddziel} for $n=1$ but $f(1)>g(1)$. Of course, $f$ and $g$
 can not be separated by a straight line.
\end{rem}

As a consequence of Theorem~\ref{tw2} we obtain
\begin{cor}\cite[Theorem 1]{NikWas95}
 Let $f,g:I\to\R$. The following conditions are equivalent:
 \begin{enumerate}[\upshape (i)]
  \item
   there exists an affine function $h:I\to\R$ such that
   $f(x)\le h(x)\le g(x)$, $x\in I$;
  \item
   for every $x,y\in I$, $t\in\left[ 0,1\right]$
   \begin{align*}
    f\bigl(tx+(1-t)y\bigr)&\le tg(x)+(1-t)g(y)\\
    \intertext{and}
    g\bigl(tx+(1-t)y\bigr)&\ge tf(x)+(1-t)f(y).
   \end{align*}
 \end{enumerate}
\end{cor}

\begin{proof}
The above ineualities are equivalent to \itemref{tw2:ii} in Theorem~\ref{tw2} (cf. the proof of Corollary~\ref{wn1}).
\end{proof}

\section{Applications}

One can verify that $f:I\to\R$ is convex iff for every
$x_{0},x_{1},x_{2}\in I$ such that $x_{0}<x_{1}<x_{2}$
\[f(x_{0})\ge c_{1}(x_{0};x_{1},x_{2})f(x_{1})+
               c_{2}(x_{0};x_{1},x_{2})f(x_{2}).\]
We adopt the following definition (cf. \cite[\S 83]{RobVar73},
\cite{Cie59}, \cite{Pop44}, \cite{Pop34}).

\begin{defn*}
Let $n\in\N$. The function $f:I\to\R$ is
$n$-\emph{convex} iff for every
$x_{0},x_{1},\dots,x_{n+1}\in I$ such that
$x_{0}<x_{1}<\dots<x_{n+1}$
\[ (-1)^{n}f(x_{0})\le (-1)^{n}\sum_{i=1}^{n+1}
   c_{i}(x_{0};x_{1},\dots,x_{n+1})f(x_{i}).\]
$f$ is $n$-\emph{concave} iff $(-f)$ is $n$-convex.
\end{defn*}

If $f$ is both $n$-convex and $n$-concave then $f$ is a polynomial
belonging to ${\mathscr{P}}_{n}$ (going through the points $(x_{i},f(x_{i}))$,
$i=0,1,\dots,n+1$).
\begin{cor}\label{wn3}
 Let $n\in\N$. If $f:I\to\R$ is $n$-convex, $g:I\to\R$ is $n$-concave and
 $f(x)\le g(x)$, $x\in I$, then there exists a polynomial $w\in {\mathscr{P}}
 _{n}$ such that $f(x)\le w(x)\le g(x)$, $x\in I$.
\end{cor}

\begin{proof}
Fix any $x_{0},x_{1},\dots,x_{n+1}\in I$ such that $x_{0}\le
x_{1}<\dots<x_{n+1}$. If $n$ is even then by $n$-convexity of $f$
\begin{align*}
 f(x_{0})&\le\sum_{i=1}^{n+1}c_{i}(x_{0};x_{1},\dots,x_{n+1})f(x_{i})\\
         &\le\sum_{\substack{i=1\\ i\text{\ \normalfont{odd}}}}^{n+1}
              c_{i}(x_{0};x_{1},\dots,x_{n+1})g(x_{i})
          +\sum_{\substack{i=1\\ i\text{\ \normalfont{even}}}}^{n+1}
            c_{i}(x_{0};x_{1},\dots,x_{n+1})f(x_{i}).
\end{align*}
If $n$ is odd then by $n$-concavity of $g$
\begin{align*}
 f(x_{0})&\le g(x_{0})\le\sum_{i=1}^{n+1}c_{i}(x_{0};x_{1},\dots,x_{n+1})g(x_{i})\\
         &\le\sum_{\substack{i=1\\ i\text{\ \normalfont{odd}}}}^{n+1}
              c_{i}(x_{0};x_{1},\dots,x_{n+1})g(x_{i})
            +\sum_{\substack{i=1\\ i\text{\ \normalfont{even}}}}^{n+1}
              c_{i}(x_{0};x_{1},\dots,x_{n+1})f(x_{i}).
\end{align*}
The proof of the second ineuality in \eqref{oddziel} is analogous.
Theorem~\ref{tw2} completes the proof.
\end{proof}

In the same way we get
\begin{cor}\label{wn4}
 Let $n\in\N$. If $f:I\to\R$ is $n$-concave, $g:I\to\R$ is $n$-convex and
 $f(x)\le g(x)$, $x\in I$, then there exists a polynomial $w\in {\mathscr{P}}
 _{n}$ such that $f(x)\le w(x)\le g(x)$, $x\in I$.
\end{cor}

For $n=1$ the above two results are well known and
Corollaries~\ref{wn3}~and~\ref{wn4} generalize them.

Finally we prove a stability result for polynomials (cf. a
Hyers-Ulam stability theorem for affine functions in \cite{NikWas95}). First
observe that if $n\in\N$ and $w(x)=1$, $x\in I$, then
$w\in {\mathscr{P}}_{n}$ and for every  different points
$x_{1},\dots,x_{n+1}\in I$ \eqref{Lagrange} has the form
\[\sum_{i=1}^{n+1} c_{i}(x;x_{1},\dots,x_{n+1})=1,\;\;x\in I.\]
\begin{cor}
 Let $n\in\N$, $\varepsilon >0$ and $f:I\to\R$. If for every
 $x_{0},x_{1},\dots,$ \\ $x_{n+1}\in I$ such that
 $x_{0}\le x_{1}<\dots<x_{n+1}$
  \begin{equation}\label{modul1}
    \biggl| f(x_{0})-\sum_{i=1}^{n+1} c_{i}(x_{0};x_{1},\dots,x_{n+1})
    f(x_{i})\biggr| \le\varepsilon
  \end{equation}
 then there exists a polynomial $w\in {\mathscr{P}}_{n}$ such that
 \begin{equation}\label{modul2}
  \left| f(x)-w(x)\right| \le\frac{\varepsilon}{2},\;\;x\in I.
 \end{equation}
\end{cor}

\begin{proof}
If $f$ satisfies \eqref{modul1} then $(ii)$ in Theorem~\ref{tw2}
holds for $g(x)=f(x)+\varepsilon$, $x\in I$. So there exists a polynomial
$\varphi\in {\mathscr{P}}_{n}$ such that $f(x)\le\varphi (x)\le f(x)+\varepsilon$,
$x\in I$. For \[w(x)=\varphi (x)-\frac{\varepsilon}{2},\;\;x\in I\]
we obtain \eqref{modul2}.
\end{proof}

\bibliographystyle{amsplain}
\bibliography{was_own,was_pub}

\end{document}